\documentclass[12pt,leqno]{amsart}%
\usepackage{amsmath}
\usepackage[
margin=1.8in,
]{geometry}
\usepackage{graphicx}
\usepackage{amsfonts}
\usepackage{amssymb}%
\usepackage[utf8]{inputenc}
\usepackage{comment}
\usepackage{hyperref}
\usepackage{mathtools}
\usepackage{ dsfont }

\newtheorem{theorem}{Theorem}[section]
\theoremstyle{plain}

\newtheorem{corollary}{Corollary}[section]

\numberwithin{equation}{section}
\theoremstyle{definition}

\newtheorem{remark}{Remark}[section]

\def\pd{\partial}
\def\re{\mathbb{R}}

\newcommand{\eqal}[1]{\begin{equation}\begin{aligned}#1\end{aligned}\end{equation}}

\begin{document}
\title[Lagrangian mean curvature equation]{Gradient estimates for the Lagrangian mean curvature equation with critical and supercritical phase}
\author{Arunima Bhattacharya}
\address{Department of Mathematics\\
University of Washington, Seattle, WA, U.S.A.}
\email{arunimab@uw.edu}

\author{Connor Mooney}
\address{Department of Mathematics\\
University of California, Irvine, CA, U.S.A.}
\email{mooneycr@math.uci.edu}

\author{Ravi Shankar}
\address{Department of Mathematics\\
Princeton University, Princeton, NJ, U.S.A.}
\email{rs1838@princeton.edu}

\begin{abstract}
  In this paper, we prove interior gradient estimates for the Lagrangian mean curvature equation, if the Lagrangian phase is critical and supercritical and $C^{2}$.  Combined with the a priori interior Hessian estimates proved in \cite{ABhess, AB2d}, this 
  solves the Dirichlet boundary value problem for the critical and supercritical Lagrangian mean curvature equation with $C^0$ boundary data.  We also provide a uniform gradient estimate for lower regularity phases that satisfy certain additional hypotheses.

\end{abstract}


\maketitle

\section{ Introduction}

In this paper, we study a priori interior gradient estimates in all dimensions for the Lagrangian mean curvature equation
\begin{equation}
\label{s}
    F(D^{2}u)=\sum _{i=1}^{n}\arctan \lambda_{i}=\psi(x), \qquad x\in B_1(0)\subset\re^n,
\end{equation}
under the assumption that $|\psi|\geq (n-2)\frac{\pi}{2}$.  Here, $u:B_1\to\re$ has gradient $Du$ and Hessian matrix $D^2u$, with eigenvalues $\lambda_i$.  We will denote $B_r=B_r(0)$ throughout.

\smallskip
When the phase $\psi$ is constant, denoted by $c$, $u$ solves the special Lagrangian equation 
 \begin{equation}
 \label{s1}
\sum _{i=1}^{n}\arctan \lambda_{i}=c,
\end{equation} 
or equivalently,
\[ 
\cos c \sum_{1\leq 2k+1\leq n} (-1)^k\sigma_{2k+1}-\sin c \sum_{0\leq 2k\leq n} (-1)^k\sigma_{2k}=0.
\]
Equation \eqref{s1} originates in the special Lagrangian geometry by Harvey-Lawson \cite{HL}.  The Lagrangian graph $(x,Du(x)) \subset \mathbb{R}^n\times\mathbb{R}^n$ is called special
when the argument of the complex number $(1+i\lambda_1)...(1+i\lambda_n)$,
or the phase $\psi$, is constant, and it is special if and only if $(x,Du(x))$ is a
(volume minimizing) minimal surface in $(\mathbb{R}^n\times\mathbb{R}^n,dx^2+dy^2)$ \cite{HL}.  

\smallskip
More generally, for \eqref{s}, it was shown in \cite[(2.19)]{HL} that the mean curvature vector $\vec H$ of the Lagrangian graph $(x,Du(x))$ is $J\nabla_g\psi$, where $\nabla_g$ is the gradient, and $J$ is the almost complex structure on $\re^n\times\re^n$.  Note that $|\vec H|_g$ is bounded for $\psi\in C^1$.  In the complex setting, a local version of the deformed Hermitian-Yang-Mills equation for a holomorphic line bundle over a compact K\"ahler manifold is represented by equation (\ref{s}).

\smallskip
The notions of critical and supercritical phases were introduced by Yuan \cite{YGlobal}.  The Lagrangian angle $\theta(\lambda)=\sum_i\arctan\lambda_i$ is critical if $|\theta|=(n-2)\pi/2$ and supercritical if $|\theta|>(n-2)\pi/2$.  We recall that the variable phase $\psi(x)$ is called critical and supercritical if $|\psi(x)|\ge(n-2)\pi/2$, and supercritical if $|\psi(x)|\ge(n-2)\pi/2+\delta$ for some $\delta>0$.   It was shown in \cite[Lemma 2.1]{YGlobal} that the level sets $\{\lambda:\theta=c\}$ are convex for critical and supercritical phases.  In particular, there are Evans \cite{EK}-Krylov \cite{KrPara}-Safonov\cite{Saf84,SafC2a} $C^{2,\alpha}$ estimates if $D^2u$ is bounded, and $\psi(x)$ is H\"older continuous.

\smallskip
In this paper, for $C^2$ critical and supercritical phases, we solve the Dirichlet problem for $C^0$ boundary data by establishing the missing interior gradient estimates. Interior Hessian estimates for supercritical $C^{1,1}$ phases were shown by Bhattacharya in \cite[Theorem 1.1]{ABhess}; interior Hessian estimates for critical and supercritical phases follow verbatim from the calculations done in \cite{ABhess} (see \cite[Remark 2.1]{AB2d}); interior gradient estimates for supercritical $C^{1}$ phases were derived in \cite[Theorem 1.2]{ABhess}. 

\medskip
Our main result is a gradient estimate for arbitrary $C^2$ critical and supercritical phases.

\begin{theorem}\label{thm:gradC2}
Let $u$ be a $C^3(\overline{B_1})$ solution of (\ref{s}) on $B_{1}(0)\subset \mathbb{R}^{n}$, where $\psi\in C^{2}(B_1)$ satisfies $\psi\geq (n-2)\frac{\pi}{2}$.  Then
\begin{equation}
\label{est:qua}
    |Du(0)|\leq C(n,||D^2\psi||_{L^\infty(B_{1})})\left(1+(\text{osc}_{B_1}u)^2\right).
\end{equation}
\end{theorem}

We state the following Hessian estimate combining \cite[Theorem 1.1]{ABhess} and \cite[Remark 2.1]{AB2d}.
\begin{theorem}\label{thm:hess}
Let $u$ be a $C^4$ solution of (\ref{s}) on $B_{R}(0)\subset \mathbb{R}^{n}$, where $\psi\in C^{2}(B_{R}) $, and $\psi\geq (n-2)\frac{\pi}{2}$. Then we have 
\begin{equation}
    |D^2u(0)|\leq C\exp [C\max_{B_R(0)}|Du|^{2n-2}/R^{2n-2}] \label{H}
\end{equation}
where $C$ is a positive constant depending on $||\psi||_{C^{2}(B_{R})}$ and $n$.
\end{theorem}

As an application, we solve the following Dirichlet boundary value problem with $C^{0}$ boundary data.  

\begin{corollary}\label{thm:dbvp}
Suppose that $\phi\in C^{0}(\partial \Omega)$ and $\psi: \overline \Omega\rightarrow [(n-2)\frac{\pi}{2}, n\frac{\pi}{2})$ is in $C^{2}(\overline \Omega)$, where $\Omega$ is a uniformly convex, bounded domain in $\mathbb{R}^{n}$. Then there exists a unique solution $u\in C^3(\Omega)\cap  C^{0}(\overline\Omega)$ to the Dirichlet  problem
\begin{align}
		\begin{cases}
		F(D^{2}u)=\sum _{i=1}^{n}\arctan \lambda_{i}=\psi(x) \text{ in }  \Omega\\
		u=\phi \text{ on }
		\partial \Omega
		\end{cases}\label{lab}
		\end{align}
\end{corollary}
\noindent The solution $u$ is, in fact, in $C^{3,\alpha}(\Omega)$ for any $\alpha \in (0,\,1)$, by classical uniformly elliptic theory.

\smallskip

The Dirichlet problem for a broad class of fully nonlinear, elliptic equations of the form $ F(\lambda[D^2u])=f(x)$ was first studied by Caffarelli-Nirenberg-Spruck in \cite{CNS}, where they proved the existence of
classical solutions under various hypotheses on the function $F$ and the domain $\Omega$. In \cite{HL09}, Harvey-Lawson studied the Dirichlet problem for fully nonlinear, degenerate elliptic equations of the form $F(D^2u)=0$ on a smooth bounded domain in $\mathbb{R}^n$. The existence and uniqueness of continuous viscosity solutions to the Dirichlet problem for (\ref{s1}) with continuous boundary data was shown in \cite{HL09,YnotesSummer08}; see also \cite{ABdir}. In \cite{BrW}, Brendle-Warren studied a second
boundary value problem for the special Lagrangian equation.
\smallskip

For subcritical phases $|c|<(n-2)\pi/2$, interior regularity is not understood.  For critical $|c|=(n-2)\pi/2$ and supercritical $|c|>(n-2)\pi/2$ phases, interior gradient estimates were established by Warren-Yuan \cite{WY2d,WYsuper}, and also Yuan's unpublished notes from 2015.  Interior Hessian estimates for dimension $n=2$ were shown by Heinz, for $|c|=\pi/2$ in dimension $n=3$ by Warren-Yuan \cite{WY093d}, and for general dimension $|c|\ge (n-2)\pi/2$ by Wang-Yuan \cite{WdYhess}; see also Li \cite{Lcomp} for a compactness approach and Zhou \cite{ZhouHess} for estimates requiring Hessian constraints which generalize criticality.  Because the level set of the PDE is convex for critical and supercritical phases, the Evans-Krylov theory yields interior analyticity.  The singular $C^{1,\alpha}$ subcritical phase solutions by Nadirashvili-Vl\u{a}du\c{t} \cite{NVsing} and Wang-Yuan \cite{WdYsing} show that interior regularity is not possible for subcritical phases, without an additional convexity condition, as in Bao-Chen \cite{BCconvex}, Chen-Warren-Yuan \cite{CWY}, and Chen-Shankar-Yuan \cite{CSY}, and that the Dirichlet problem is not classically solvable for arbitrary smooth boundary data.  
  Interior gradient estimates for continuous boundary data are widely open. Global gradient estimates requiring Lipschitz boundary data were shown by \cite{Siyuan}. Homogeneous viscosity solutions of degree less than two were shown to not exist by Nadirashvili-Yuan \cite{NYhom}. 
The non-existence result of Mooney \cite{Moon} shows that counterexamples for interior $C^1$ regularity may be difficult to construct.
\smallskip

If the Lagrangian angle is not necessarily constant, then less is understood.  
In \cite{HL19}, Harvey-Lawson introduced a condition called “tameness” on the operator $F$, which is a little stronger than strict ellipticity and allows one to prove comparison.  In \cite{HL21}, tamability was established for the supercritical Lagrangian mean curvature equation. In \cite{CP}, Cirant-Payne established comparison principles for the Lagrangian mean curvature equation provided the Lagrangian phase is restricted to the intervals $((n-2k)\frac{\pi}{2},(n-2(k-1))\frac{\pi}{2})$ where $1\leq k\leq n$, which in turn solves the Dirichlet problem on these intervals as shown in \cite[Theorem 6.2]{HL21}.
Hessian estimates for convex smooth solutions with $C^{1,1}$ phase $\psi=\psi(x)$ were obtained by Warren \cite[Theorem 8]{WTh}.  For convex viscosity solutions, interior regularity was established for $C^2$ phases; see Bhattacharya-Shankar \cite{BS1,BS2}. For supercritical phases $|\psi(x)|\ge(n-2)\pi/2+\delta$, there is a comparison principle, and the Dirichlet problem was solved in Collins-Picard-Wu \cite{CPW}, Dinew-Do-T{\^o} \cite{DDT}, Bhattacharya \cite{ABdir}, and interior gradient estimates were established in \cite{ABhess}. Interior Hessian estimates for supercritical phases were established in \cite{ABhess}.
Interior Hessian estimates for critical and supercritical phases $|\psi(x)|\geq (n-2)\frac{\pi}{2}$, follow verbatim from the calculations done in \cite{ABhess} (see \cite[Remark 2.1]{AB2d}): The proof of the Hessian estimate in \cite[Theorem 1.1]{ABhess} does not require a negative lower bound on the lowest eigenvalue. For supercritical phases in dimension $n=2$, a simplified proof \cite{AB2d} was given for interior Hessian estimates using the super-isoperimetric inequality of Warren-Yuan \cite{WY2d}, avoiding the Michael-Simon mean value inequality \cite{MS}. In the case that $\phi$ is Lipschitz, Corollary \ref{thm:dbvp} can be obtained by proving a global gradient estimate, as in \cite{Siyuan}. The existence of 
interior gradient estimates for the challenging borderline case of critical and supercritical phase has until now remained open. In this paper, we successfully solve this problem for $C^2$ phases.
\smallskip

Our approach to prove interior gradient estimate Theorem \ref{thm:gradC2} accounts for the smallness of the gradient of the phase near its minimizing, critical values, using a pointwise interpolation inequality \cite[Equation (3.11), pg. 19]{NT70}, see also \cite[Lemma 7.7.2]{Ho1}, valid for $C^2$ phases.  
%
%
%
%
%
%
%
%
%
%
%
%
%
%
For constant phases, the gradient estimate is established using a maximum principle inequality.  The variable phase contribution to the inequality is a ``bad term" depending on the phase's gradient.  
Although the 
PDE's ellipticity degenerates at the critical phase, making the bad term large, the smallness of the gradient at such points provides a balance.  Our proof, more generally, shows that an interior gradient estimate holds when $\psi$ satisfies a certain first order differential inequality; see Remark \ref{rem:general}.  Such an inequality is valid when the phase is any of $C^2$, semi-concave, concave, or a supersolution of the infinity-Laplace equation; see Remark \ref{rem:example}.
\smallskip

On the other hand, the gradient vanishes at slower rates for $C^{1,\alpha}$ phases, and does not appear to balance the degeneration of the ellipticity in our proof.  
But we note that certain H\"older continuous phases allow for gradient estimates; see Remark \ref{rem:xa}.  In such cases, the phase separates from the critical value at such a large speed 
that the solution is nearly semi-convex, as in the supercritical case.

\medskip
\textbf{Acknowledgements.} A. Bhattacharya is grateful to Yu Yuan for bringing attention to this problem and for many helpful discussions. R. Shankar thanks Yu Yuan for helpful suggestions about the paper.  A. Bhattacharya would like to thank Pengfei Guan for pointing out that the interpolation inequality appears in \cite{NT70}, as well as Tristan Collins for helpful conversations.
C. Mooney gratefully acknowledges the support of NSF grant DMS-1854788, NSF CAREER grant DMS-2143668, an Alfred P. Sloan Research Fellowship, and a UC Irvine Chancellor’s Fellowship.

\section{Proof of the gradient estimate}

We modify the pointwise proof of \cite{WY2d} and \cite{ABhess} to bridge the constant critical phase estimate of \cite{WY2d} with the supercritical estimate of \cite{ABhess}.  The difference comes from how to treat the bad term involving $D\psi$.

\medskip
\textbf{Notation:}  we denote $a\sim b$ if $ca\le b\le Ca$, and $a\lesssim b$ if $a\le Cb$.  Here, $c$ and $C$ are positive constants depending on $n$.  We will denote $a\lesssim_\psi b$ if the constant also depends on $\|D^2\psi\|_{L^\infty(B_{1})}$.  We denote $a\ll_n 1$ if we are choosing a small fixed constant $a$ depending on $n$.  We assume summation under repeated indices unless otherwise indicated. 

\medskip
Let $M=\text{osc}_{B_1}u>0$; replacing $u$ with $u-\min_{B_1}u+M$, we assume that 
\eqal{
M\le u\le 2M\qquad \text{in }B_1.
}
Let $w=\eta|Du|+Au^2/2$, where $\eta=1-|x|^2$, and $A=3\sqrt{n}/M$.  Let $x_0\in B_1$ be where $w$ is maximized.  After a rotation, we assume that $D^2u$ is diagonal, with $u_{ii}=\lambda_i$.  Let us assume that $u_n\ge |Du|/\sqrt{n}>0$.  Then for each $k$, at the max point $x_0$,
\eqal{
\label{Dw}
0=\pd_kw(x_0)=\eta\frac{u_k\lambda_k}{|Du|}+\eta_k|Du|+Auu_k.
}
Since $A=3\sqrt{n}/M$ is sufficiently large, it follows that
\eqal{
\eta\lambda_n\frac{u_n}{|Du|}\in -(c(n),C(n))|Du|.
}
It follows that $\lambda_n<0$ and
\eqal{
\label{eta}
\eta\sim\frac{|Du|}{|\lambda_n|}.
}
Since $|Du|\lesssim \eta|\lambda_n|$, we may assume that $|\lambda_n|> 1$, since otherwise the estimate is done.  
Moreover, as shown in \cite[Lemma 2.5]{WdYhess}, we know that $\lambda_k\ge |\lambda_n|$ for $k<n$ follows from $\psi\ge (n-2)\pi/2$.

\medskip
We now proceed to the second derivatives of $w$.  
Let $g=I+(D^2u)^2$ be the induced metric $dx^2+dy^2$ on $(x,Du(x))$, with $g^{-1}=(g^{ij})$ its inverse, and $g^{ij}=(1+\lambda_i^2)^{-1}\delta_{ij}$ at $x_0$.  Then at $x_0$,
\eqal{
\label{D2w}
0&\ge g^{ij}\pd_{ij}w(x_0)\\
&=\underbrace{g^{ij}\eta_{ij}|Du|}_{(C1)}+\underbrace{2g^{ij}\eta_i|Du|_j}_{(C2)}+\underbrace{\eta g^{ij}|Du|_{ij}}_{(B)}+\underbrace{Aug^{ij}u_{ij}}_{(C3)}+\underbrace{Ag^{ij}u_iu_j}_{(G)}.
}
The good term (G) absorbs the constant-phase terms (C1), (C2), and (C3) in the constant phase and supercritical cases, as in \cite{WY2d} and \cite{ABhess}.  The bad term (B) contains variable phase contributions and will require closer examination.

\medskip
The good term (G):
\eqal{
\label{(G)}
Ag^{ij}u_iu_j\ge A\frac{u_n^2}{1+\lambda_n^2}\sim A\frac{|Du|^2}{\lambda_n^2} \sim \frac{A\eta|Du|}{|\lambda_n|}.
}
The first constant phase term (C1):
\eqal{
g^{ij}\eta_{ij}|Du|=-2\sum_i\frac{1}{1+\lambda_i^2}|Du|\gtrsim -\frac{1}{\lambda_n^2}|Du|\sim -\frac{\eta}{|\lambda_n|}.
}
The second constant phase term (C2):
\eqal{
2g^{ij}\eta_i|Du|_j\stackrel{x_0}{\gtrsim}-\sum_i\frac{1}{1+\lambda_i^2}\frac{u_i\lambda_i}{|Du|}\gtrsim -\frac{1}{|\lambda_n|}.
}
The third constant phase term (C3):
\eqal{
Aug^{ij}u_{ij}=Au\sum_i\frac{\lambda_i}{1+\lambda_i^2}\gtrsim -\frac{1}{|\lambda_n|}. 
}
The bad term (B), using third derivative calculation \cite[(2.4)]{ABhess}:
\eqal{
\label{(B)}
\eta g^{ij}|Du|_{ij}&=\eta g^{ij}\frac{u_{ijk}u_k}{|Du|}+\eta \sum_ig^{ii}\frac{(|Du|^2-u_i^2)\lambda_i^2}{|Du|^3}\\
&\ge \eta \psi_k\frac{u_k}{|Du|}\gtrsim -\eta|D\psi|.
}
We thus need to bound this inequality at $x_0$:
\eqal{
\label{etaDu}
\eta|Du|\le C(n) M(1+\eta\,|D\psi|\,|\lambda_n|).
}


\smallskip
\noindent
Letting $\phi=\psi-(n-2)\pi/2\ge 0$, we apply the pointwise interpolation inequality for nonnegative $C^2$ functions in \cite[Lemma 7.7.2]{Ho1} on $B_\delta(x_0)$, where $\delta=1-|x_0|$;
\eqal{
\label{Dtheta}
|D\phi(x_0)|^2&\le \frac{\phi(x_0)^2}{(1-|x_0|)^2}+2\|D^2\phi\|_{L^\infty(B_1)}\phi(x_0)\\
&\lesssim_\psi \frac{\phi(x_0)^2}{\eta^2}+\phi(x_0),
}
where $\|D^2\phi\|_{L^\infty(B_1)}$ denotes the maximum of the absolute values of the eigenvalues of $D^2\phi$.  Let us now recall the following algebraic inequality, valid for $\lambda_n<0$ and $\lambda_k>0$ for $k<n$:
\eqal{
\label{alg}
\psi&=(n-1)\frac{\pi}{2}-\sum_{i<n}\arctan(\frac{1}{\lambda_i})-\frac{\pi}{2}+\arctan(-\frac{1}{\lambda_n})\\
&\le (n-2)\frac{\pi}{2}+\frac{1}{|\lambda_n|}.
}
Substituting this information into \eqref{etaDu}, combined with \eqref{eta}, yields
\eqal{
\label{etaDusub}
M^{-1}\eta |Du|&\lesssim_\psi 1+\eta|\lambda_n|^{1/2}\lesssim 1+(\eta |Du|)^{1/2}.
}
It follows that 
\eqal{
\eta|Du|\lesssim_\psi M+M^2\lesssim C(n,\|\psi\|_{C^2(B_{1.5})})(1+(\text{osc}_{B_1}u)^2).
}
The $Au^2/2$ term in $w(x_0)$ and the estimate of $w$ on $\pd B_1$ are subordinate to this estimate, so we conclude the proof.
$\qed$

\begin{remark}
\label{rem:refine}
It is straightforward to refine the $|\lambda_n|<1$ case and thereby improve the estimate \eqref{est:qua} 
to the following
\eqal{
&|Du(0)|\le C(n,\|D^2\psi\|_{L^\infty(B_1)})(\text{osc}_{B_1}u+(\text{osc}_{B_1}u)^2).
}
\end{remark}

\begin{remark}
\label{rem:general}
More generally, let $\psi\in C^1(B_1)$ be critical and supercritical, or $\psi-(n-2)\pi/2=:\phi\ge 0$, and also satisfy the following first order differential inequality on $B_1$:
\eqal{
\label{DthetaGeneral}
|D\phi|\le \eta f(\eta^{-2}\phi),
}
where $f(t)\searrow 0$ as $t\searrow 0$, and $\eta=1-|x|^2$.  Then a $C^1$ estimate is valid for $C^3(\overline{B_1})$ solutions of \eqref{s}:
\eqal{
\label{est:gen}
|Du(0)|\le C(n,f,\text{osc}_{B_1}u).
}
To prove the estimate, we insert \eqref{DthetaGeneral} in the determining inequality \eqref{etaDu} and use \eqref{alg} and \eqref{eta}.  We obtain at $x_0$, using that $f$ is increasing,
\eqal{
M^{-1}\,\eta|Du|&\lesssim 1+\eta^2|\lambda_n|f\left(\frac{1}{\eta^2|\lambda_n|}\right)\\
&\lesssim 1+\eta |Du| f\left(\frac{C(n)}{\eta|Du|}\right).
}
If $\eta|Du|\ge C(n,f,M)=:H$ for large enough $H$ such that $f(C(n)/H)\ll_{n} M^{-1}$, then $M^{-1}\eta|Du|\lesssim 1$, and the estimate follows.  In the alternative case that $\eta|Du|\le C(n,f,M)$, the estimate \eqref{est:gen} is already done.
\end{remark}

\begin{remark}
\label{rem:example}
Let us list some examples of phases which satisfy a first order inequality of the form \eqref{DthetaGeneral}.  

\smallskip
\noindent
1. For $f(t)^2=t^2+Ct$, we recover the $C^2$ interpolation inequality \eqref{Dtheta}.  Note that general $C^{1,\alpha}$ phases fail to satisfy the inequality \eqref{DthetaGeneral}.

\smallskip
\noindent
2. Interpolation inequality \eqref{Dtheta} can be generalized to phases $\psi\in C^1(B_1)$ which are semi-concave, with $D^2\psi\le KI$ for some $K> 0$.  In this case, the dependence on $\|D^2\psi\|_{L^\infty(B_1)}$ is replaced with $K$.  Indeed, by semi-concavity, there holds for $x_0\in B_1(0)$ and $x\in B_\delta(x_0)$:
\eqal{
\label{semi}
0\le \phi(x_0)+(x-x_0)\cdot D\phi(x_0)+K|x-x_0|^2/2.
}
The proof in \cite[Lemma 7.7.2]{Ho1} can then be repeated verbatim.  This generalizes Theorem \ref{thm:gradC2} to semi-concave phases.

\smallskip
\noindent
3. The choice $f(t)=2t$ corresponds to $\psi\in C^1(B_1)$ concave.  Choosing $x-x_0=-(1-|x_0|)D\phi(x_0)/|D\phi(x_0)|$ with $K=0$ in \eqref{semi} gives
\eqal{
\label{Dthetaconcave}
|D\phi(x_0)|\le \frac{\phi(x_0)}{1-|x_0|}\le \frac{2\phi(x_0)}{\eta}.
}
This is the first term in \eqref{Dtheta}, so as in \eqref{etaDusub}, we obtain $\eta|Du|\le C(n)M$.  We thus obtain the linear estimate
\eqal{
\label{est:lin}
|Du(0)|\le C(n)(1+\text{osc}_{B_1}u).
}
This can be improved to $|Du(0)|\le C(n)\text{osc}_{B_1}u$, as in Remark \ref{rem:refine}.  One novelty here is the independence of $\psi$.  For example, if 
\eqal{
\psi(x)=(n-2)\frac{\pi}{2}+\epsilon(1-|x|^{1+\alpha})
}
for some $\epsilon,\alpha\in(0,1)$, then \eqref{est:lin} is independent of $\epsilon$.  The interior gradient estimate for $C^1$ supercritical phases in \cite{ABhess} would degenerate as $\epsilon\to 0$.

\smallskip
\noindent
4. Suppose that $\phi(x)\ge 0$ is a $C^1(B_1)$ viscosity supersolution of the infinity-Laplace equation/Aronsson's equation:
\eqal{
D^2\phi(D\phi,D\phi)=\phi_{ij}\phi_i\phi_j\le 0.
}
Then using comparison with cones, there is a pointwise estimate \cite[Lemma 2.5]{CEG} for the gradient:
\eqal{
\label{Daron}
|D\phi(x)|\le \frac{\phi(x)}{1-|x|}.
}
In fact, this is concavity inequality \eqref{Dthetaconcave}, and this corresponds to $f(t)=2t$ in \eqref{DthetaGeneral}.  We conclude that a linear gradient estimate \eqref{est:lin} is valid.
\end{remark}

\begin{remark}
\label{rem:xa}
If $u$ is a viscosity solution of \eqref{s} for H\"older phase 
\eqal{
\psi=(n-2)\frac{\pi}{2}+|x|^{\alpha},
}
where $0<\alpha<1$, the function
\eqal{
u(x)+C|x|^{2-\alpha}
}
is convex, if $C(\alpha)$ is large enough.  This follows from the algebraic relation \eqref{alg}, which gives $|\lambda_{min}|<|x|^{-\alpha}$.  
It follows that $u(x)$ is locally Lipschitz continuous.  
\end{remark}

\bibliographystyle{amsalpha}
\bibliography{lmc}

\end{document}